\documentstyle[12pt,leqno]{article}
\title{
 \bf Monopoles and contact 3-manifolds}
\author{
  Jih-Hsin Cheng 
\thanks{
Research supported in part by National Science Council grant NSC 
87-2115-M-001-006 (R.O.C.).}
and Hung-Lin Chiu}
\date{}

\begin{document}
\maketitle
{\normalsize \setlength{\baselineskip}{18pt}
\begin{abstract}
We propose the study of some kind of monopole equations directly 
associated with a contact structure. Through a rudimentary analysis
about the solutions, we show that a closed contact 3-manifold
with positive Tanaka-Webster curvature and vanishing torsion must be either
not symplectically semifillable or having torsion Euler class of the
contact structure.
\end{abstract}



\section{Statement of results}

In this paper we propose some kind of monopole equations directly
associated to a contact structure. By studying the solutions of
these equations, we can draw a conclusion about the underlying
contact structure.

Given an oriented contact structure $\xi$ on a closed
(compact without boundary) 3-manifold $M$, we can talk about
$spin^c$-structures on $\xi$ or ${\xi}^{\star}$. (see $\S$2
for the definition) Furthermore, associated to an oriented
pseudohermitian structure, we have the so-called
canonical $spin^c$-structure
$c_{\xi}$. With respect to $c_{\xi}$, we consider the equations
(3.9) for our ``monopole'' $\Phi$ coupled to the ``gauge field''
$A$. Here $A$, the $spin^c$-connection, is required to be
compatible with the pseudohermitian connection on $M$. The
Dirac operator $D_{\xi}$ relative to $A$ is identified with
a certain boundary $\bar {\partial}$-operator 
${\sqrt 2}({\bar {\partial}}_b^{a}+
({\bar {\partial}}_b^{a})^{\star})$. (cf.(3.10)) In terms of
components $({\alpha},{\beta})$ of $\Phi$, (3.9) is equivalent
to

\begin{flushleft}(3.11)\hspace{.8in}
$
\left\{ \begin{array}{c}
       ({\bar {\partial}}_b^{a}+
({\bar {\partial}}_b^{a})^{\star})({\alpha}+{\beta})=0\\
(or\:{\alpha}_{,{\bar 1}}^{a}=0,\:{\beta}_{{\bar 1},1}^{a}=0)\\
da(e_{1},e_{2})-{\cal W}=|{\alpha}|^{2}-|{\beta}_{\bar 1}|^{2}
\end{array} \right.
$
\end{flushleft}

On the other hand, there are notions of symplectic fillability
and symplectic semifillability in the study of contact
structures due to Eliashberg. (e.g.,[ET],[Kro])
A contact 3-manifold $(M,{\xi})$ is
symplectically fillable if $\xi$ is positive (i.e. ${\theta}{\wedge}
d{\theta}>0$ for any contact form $\theta$) with respect to the induced
orientation on $M$ as the boundary of the canonically oriented
symplectic 4-manifold $(X,{\Omega})$ and ${\Omega}|_{\xi}{\neq}0$.
If $M$ consists of a union of
components of such a boundary, then it is symplectically
semifillable. 

Let $e({\xi})$ denote the Euler class of the contact bundle 
$\xi$. We say the equations (3.11) have nontrivial solutions
if $\alpha$ and $\beta$ are not identically zero simultaneously.
Our first step to understand
the equations (3.11) is the following result.

\medskip

$\bf {Theorem\:A}$. Suppose there is an oriented 
pseudohermitian structure
with vanishing torsion on a closed 3-manifold $M$ with an
oriented contact structure $\xi$.
Also suppose $\xi$ is symplectically semifillable, and
$e({\xi})$ is not a torsion class. Then 
the equations (3.11) (for the canonical $spin^c$-structure
$c_{\xi}$) have nontrivial solutions.

\medskip

We remark that our $M$ in Theorem A must be a Seifert fibre
space with even first Betti number by an argument of
Weinstein. ([CH])
The idea of proving Theorem A goes as follows. The contact structure
$\xi$ being symplectically semifillable implies that its
Euler class $e({\xi})$
is a so-called monopole class in Kronheimer's terminology. 
(see Corollary 5.7 in [Kro]) That is to say, $e({\xi})$
arises as the first Chern class of a usual (i.e. on $TM$ or 
$T^{\star}M$) $spin^c$-structure for which the usual
Seiberg-Witten equations admit a solution for every choice
of Riemannian metric on $M$. By choosing a suitable family
of Riemannian metrics adapted to our pseudohermitian 
structure, we prove that the associated solutions admit
a subsequence converging to a nontrivial solution of our
equations (3.11). (see $\S$4 for details)

On the other hand, associated to an oriented
pseudohermitian structure on a contact manifold is the
notion of the so-called Tanaka-Webster curvature $\cal W$.
([Tan],[Web],[CL], see also $\S$5) The Weitzenbock-type formula
tells a nonexistence result: (see $\S$3 for details)

\medskip

$\bf {Theorem\:B}$. Let $(M,{\xi})$ be a closed 3-manifold
with an oriented contact structure $\xi$. Suppose there is
an oriented pseudohermitian structure on
$(M,{\xi})$ with ${\cal W}>0$. Then
the equations (3.11) have no nontrivial solutions with

\begin{equation}
{\alpha}_{,0}^{a}=0,\:{\beta}_{{\bar 1},0}^{a}=0.
\end{equation}

\medskip

The solution we find for Theorem A actually satisfies the
condition (1.1). Therefore by Theorems A and B, we can
conclude 

\medskip

$\bf {Corollary\:C}$. Let $(M,{\xi})$ be a closed 3-manifold
with an oriented contact structure $\xi$. Suppose there is
an oriented pseudohermitian structure on
$(M,{\xi})$ with vanishing torsion and ${\cal W}>0$.
Then either $\xi$ is not symplectically semifillable or $e({\xi})$ is a
torsion class.

\medskip

We remark that Rumin ([Rum]) proved that $M$
must be a rational homology sphere under the conditions in
Corollary C by a different method. On the other hand we feel
that we haven't made use of the full power of equations (3.11).
Also note that Eliashberg gives a complete
list of classes in $H^{2}(L(p,1),Z)$, which can be realized
as Euler classes of fillable contact structures on the lens
spaces $L(p,1)$. ([Eli]) 

During the preparation of this paper we noticed that
Nicolaescu had a similar consideration of the so-called adiabatic
limit as in our proof of Theorem A. ([Nic]) But our viewpoint
is sufficiently different from his. Also we noticed that Kronheimer
and Mrowka ([KM],[Kro]) had studied contact
structures on 3-manifolds via 4-dimensional monopole invariants
introduced by Seiberg and Witten. ([Wit])

Since our Dirac operator $D_{\xi}$ (also $da(e_{1},e_{2})$)
is not elliptic (not even subelliptic) from our knowledge
about ${\bar {\partial}}_b$-operator, we do not know how to
deal with the solution space of (3.11) in general.
\medskip

$\bf {Acknowledgments}$. This work was being done during the
first author's visit at Harvard University in the 97-98 
academic year. He would therefore like to thank the members
of the Mathematics Department, and especially Professor
Shing-Tung Yau, for their hospitality during his stay.

\bigskip     
  
\section{$Spin^{c}$-structures on contact bundles}
\setcounter{equation}{0}

Let $(M,{\xi})$ be a smooth contact 3-manifold with oriented contact
bundle $\xi$. Choose an oriented pseudohermitian structure
$(J,{\theta})$ 
compatible with $\xi$ (see $\S$5, the Appendix) so that
$h(u,v)=\frac{1}{2}d{\theta}(u,Jv)$
for $u,v{\in}{\xi}$ defines a Riemannian structure on $\xi$.
Let ${\xi}^{\star}$ denote the dual of $\xi$. The $h$ also induces
a Riemannian structure on ${\xi}^{\star}$, still denoted $h$. 
A $spin^c$-structure on $({\xi}^{\star}, h)$ (or similarly on
$({\xi}, h)$, cf.[Sal]) is a pair $(W,{\Gamma})$ where $W$
is a 2-dimensional complex Hermitian vector bundle and
${\Gamma}:{\xi}^{\star}
{\rightarrow}End(W)$ is a homomorphism which satisfies

\begin{equation}
{\Gamma}(v)^{\star}+{\Gamma}(v)=0,
{\Gamma}(v)^{\star}{\Gamma}(v)=|v|_{h}^{2}I.
\end{equation}

\noindent
Here $I$ means the identity endomorphism. Let $C^{c}({\xi}^{\star})$
denote the 
bundle of complexified Clifford algebras of ${\xi}^{\star}$.
Then $\Gamma$ extends to an algebra (bundle) isomorphism
$:C^{c}({\xi}^{\star})
{\rightarrow}End(W)$, still denoted $\Gamma$. A Hermitian connection
$\nabla$ on $W$ is called a $spin^c$-connection if there is a connection
on ${\xi}^{\star}$, also denoted $\nabla$, such that 

\begin{equation}
{\nabla}_{v}({\Gamma}(w){\Phi})={\Gamma}(w){\nabla}_{v}{\Phi}
+{\Gamma}({\nabla}_{v}w){\Phi}
\end{equation}

\noindent
for ${\Phi}{\in}C^{\infty}(M,W)$ and $w{\in}C^{\infty}(M,{\xi}^{\star}),
v{\in}C^{\infty}(M,TM)$. A $spin^c$-connection $\nabla$ on $W$ is
said to be compatible with the pseudohermitian connection
on ${\xi}^{\star}$ if it satisfies
(2.2) with ${\nabla}_{v}w$ denoting the pseudohermitian connection
induced on ${\xi}^{\star}$. (see $\S$5 and note that we'll often
view ${\xi}^{\star}$ as the
orthogonal complement of $\theta$ in $T^{\star}M$ with respect to the
adapted metric ${\theta}{\otimes}{\theta}+h$)  

Let $e^1$, $e^2$ be a positively oriented orthonormal basis of
${\xi}^{\star}$. Denote ${\varepsilon}=e^{2}e^{1}$. Then 
${\varepsilon}^{2}=-1$ and thus ${\Gamma}({\varepsilon})$
has eigenvalues ${\pm}i$. Let $W^{\pm}=\{ \Phi\in W:{\Gamma}({\varepsilon})
{\Phi}={\pm}i{\Phi}\}$. Then $W=W^{+}{\oplus}W^{-}$, and $dim_{C}W^{\pm}=1$.
Note that ${\Gamma}(v)$ maps $W^{\pm}$ to $W^{\mp}$, and every
$spin^{c}$-connection $\nabla$ on $W$ preserves subbundles $W^{+}$ and
$W^{-}$ respectively.

Next we'll define a canonical $spin^{c}$-structure and connection associated
to an oriented pseudohermitian structure $(J,{\theta})$ on our
contact manifold
$(M,{\xi})$. Let ${\Lambda}^{0,1}{\xi}^{\star}$ be the bundle of
complex 1-forms of type $(0,1)$. (a typical element is
${\theta}^{\bar 1}=e^{1}-ie^{2}$) 
Let ${\cal C}(={\Lambda}^{0,0})$ denote the trivial complex line 
bundle. Consider

\begin{equation} 
W_{can}={\cal C}{\oplus}{\Lambda}^{0,1}{\xi}^{\star}
\end{equation}

\noindent
with the natural Hermitian structure induced by $h$. Define
${\Gamma}_{can}:{\xi}^{\star}{\rightarrow}End(W_{can})$ by 

\begin{eqnarray}
{\Gamma}_{can}(e^{1}){\tau}=\frac{1}{\sqrt 2}{\theta}^{\bar 1}{\wedge}
{\tau}-{\sqrt 2}{\iota}(e_{1}){\tau} \\
{\Gamma}_{can}(e^{2}){\tau}=\frac{1}{\sqrt 2}i{\theta}^{\bar 1}{\wedge}
{\tau}-{\sqrt 2}{\iota}(e_{2}){\tau}
\end{eqnarray}

\noindent
where \{$e_1$, $e_2$\} in $\xi$ is a dual basis of \{$e^1$, $e^2$\}, and
$\iota$ denotes the interior product. The above definition is independent 
of the choice of bases. It is a direct verification that
$(W_{can}, {\Gamma}_{can})$ is a $spin^{c}$-structure on $({\xi}^{\star},
h)$. Also $W_{can}^{+}={\cal C}$, $W_{can}^{-}=
{\Lambda}^{0,1}{\xi}^{\star}$. We call $(W_{can}, {\Gamma}_{can})$
the canonical $spin^{c}$-structure on $({\xi}^{\star},h)$, denoted
$c_{\xi}$.

We know that the pseudohermitian connection preserves the subspaces
${\Lambda}^{0,k}{\xi}^{\star}$, hence $W_{can}$.
When it is restricted to $W_{can}$, we denote it by ${\nabla}_{can}$.

\medskip
$\bf {Proposition \: 2.1}$ ${\nabla}_{can}$ is a $spin^c$-connection on
$W_{can}$, which is compatible with the pseudohermitian connection
on ${\xi}^{\star}$.
\medskip

Proof: It is enough to verify (2.2) for $w=e^{1}, e^{2}$. Let $f$ be a
smooth section of ${\cal C}$, i.e. a smooth complex-valued function
on $M$. Let $v$ be a tangent vector of $M$. For simplicity, we use
${\Gamma}$, $\nabla$ instead of ${\Gamma}_{can}$, ${\nabla}_{can}$,
respectively. We compute by (2.5), (2.4)

$$ \begin{array}{c}
{\Gamma}(e^{1}){\nabla}_{v}f+{\Gamma}({\nabla}_{v}e^{1})f \\
=\frac{1}{\sqrt 2}df(v){\theta}^{\bar 1}
+\frac{1}{\sqrt 2}if{\omega}(v){\theta}^{\bar 1} \\
(write \:{\nabla}_{v}e^{1}={\omega}(v)e^{2} \: and \: 
{\nabla}_{v}e^{2}=-{\omega}(v)e^{1} \: where \: {\omega} \: is \:
the \: connection \: 1-form)\\
={\nabla}_{v}(\frac{1}{\sqrt 2}f{\theta}^{\bar 1})
={\nabla}_{v}({\Gamma}(e^{1})f)\:
(note \: that \: {\nabla}_{v}{\theta}^{\bar 1}=i{\omega}(v)
{\theta}^{\bar 1}) \end{array}
$$

For $\tau$ being a smooth section of ${\Lambda}^{0,1}{\xi}^{\star}$,
we compute

\begin{eqnarray}
{\nabla}_{v}({\Gamma}(e^{1}){\tau})&=&-{\sqrt 2}{\nabla}_{v}
({\iota}(e_{1}){\tau}) \: (by \: (2.4))\nonumber\\
 &=&-{\sqrt 2}{\nabla}_{v}({\tau}(e_{1}))\nonumber\\
 &=&-{\sqrt 2}(({\nabla}_{v}{\tau})(e_{1})+{\tau}({\nabla}_{v}e_{1}))
\nonumber\\
 &=&{\Gamma}(e^{1}){\nabla}_{v}{\tau}+{\Gamma}({\nabla}_{v}e^{1})
{\tau}\nonumber
\end{eqnarray}

Similarly we can verify (2.2) for $w=e^2$.

\begin{flushright}Q.E.D.
\end{flushright}

\medskip

Let $E$ be a Hermitian line bundle over $M$. Let $W=W_{can}{\otimes}E$,
${\Gamma}={\Gamma}_{can}{\otimes}id$. Then $(W,{\Gamma})$ defines a 
$spin^c$-structure on ${\xi}^{\star}$. (we often suppress the metric $h$)
Conversly, we have

\medskip

$\bf {Proposition \: 2.2}$ Any $spin^c$-structure $(W,{\Gamma})$
on ${\xi}^{\star}$ is
isomorphic to

$$
(W_{can}{\otimes}E,{\Gamma}_{can}{\otimes}id)
$$

\noindent for some Hermitian line bundle $E$.

\medskip

Proof: Define $\#:{\xi}{\rightarrow}{\xi}^{\star}$ by
${\#}(v)=h(v,{\cdot})$. Let ${\tilde{\Gamma}}={\Gamma}{\circ}{\#}$.
Since  ${\tilde{\Gamma}}(Jv)={\tilde{\Gamma}}(v){\Gamma}({\varepsilon})$
for $v{\in}{\xi}$, we have  ${\tilde{\Gamma}}(Jv){\Phi}=-i
{\tilde{\Gamma}}(v){\Phi}$ for ${\Phi}{\in}W^{-}$. So   
${\tilde{\Gamma}}({\cdot}){\Phi}$ is a section of the bundle
${\Lambda}^{0,1}{\xi}^{\star}{\otimes}W^{+}$.
Furthermore, the map given by

$$
{\Phi}{\rightarrow}-\frac{1}{\sqrt 2}{\tilde{\Gamma}}({\cdot}){\Phi}
$$

\noindent is a unitary isomorphism from $W^{-}$ onto
${\Lambda}^{0,1}{\xi}^{\star}{\otimes}W^{+}$.
 
Now choose $E=W^{+}$. (note that $W^{+}$ is a Hermitian line bundle)
It follows that $W^{+}{\simeq}{\cal C}{\otimes}W^{+}={\cal C}
{\otimes}E$ and $W^{-}{\simeq}{\Lambda}^{0,1}{\xi}^{\star}{\otimes}W^{+}
={\Lambda}^{0,1}{\xi}^{\star}{\otimes}E$. Also it is easy to
verify that ${\Gamma}{\simeq}{\Gamma}_{can}{\otimes}id$.
\begin{flushright}Q.E.D.
\end{flushright}

\medskip

We remark that if M is a homology sphere, then there exists one and only
one $spin^c$-structure on ${\xi}^{\star}$ (or $\xi$), 
which is the canonical one.

Let $C_{2}({\xi}^{\star})$ denote the subspace of $C({\xi}^{\star})$
(the real Clifford algebra of ${\xi}^{\star}$), consisting of elements
of degree 2.

\medskip

$\bf {Proposition \: 2.3}$ Given a $spin^c$-structure $(W,{\Gamma})$
on ${\xi}^{\star}$. Let ${\nabla}^1$, ${\nabla}^2$ be two
$spin^c$-connections on $W$. Then there exists a 1-form $\alpha$ with
value in $C_{2}({\xi}^{\star}){\oplus}i{\bf R}$ so that

$$
{\nabla}^1-{\nabla}^2={\Gamma}({\alpha})
$$

\noindent Conversly, if $\nabla$ is a $spin^c$-connection, so is 
${\nabla}+{\Gamma}({\alpha})$ for any 
$C_{2}({\xi}^{\star}){\oplus}i{\bf R}$-valued 1-form $\alpha$.

\medskip

The proof of Proposition 2.3 is similar to the
usual case for $spin^c$-structures on the tangent bundle. We
include a proof for the reference.

\medskip

Proof: Write ${\nabla}^1-{\nabla}^2=A$ for some $End(W)$-valued
1-form $A$ and express the difference of corresponding connections
on ${\xi}^{\star}$ by $a$, a $End({\xi}^{\star})$-valued 1-form.
Taking the difference of (2.2) for ${\nabla}^1$, ${\nabla}^2$ gives

$$
A(v){\Gamma}(w)-{\Gamma}(w)A(v)={\Gamma}(a(v)w)
$$

\noindent for $v{\in}TM$, $w{\in}{\xi}^{\star}$. Put $A(v)={\Gamma}
({\alpha}_{v})$ for some ${\alpha}_{v}{\in}C^{c}({\xi}^{\star})$.
Then the above formula says

\begin{equation}
{\alpha}_{v}w-w{\alpha}_{v}=a(v)w
\end{equation}

On the other hand, $A(v)$ is skew-Hermitian since ${\nabla}^1$
and ${\nabla}^2$ are Hermitian. It follows that

\begin{equation}
{\alpha}_{v}+{\tilde {{\alpha}_{v}}}=0
\end{equation}

\noindent where $\tilde {{\alpha}_{v}}$ denotes the involution
of ${\alpha}_{v}$. Now (2.6),(2.7) implies ${\alpha}_{v}{\in}
C_{2}({\xi}^{\star}){\oplus}i{\bf R}$. Let ${\alpha}(v)=
{\alpha}_{v}$. Then $\alpha$ is the required 1-form.
  
For the second part of the Proposition, we define an
$End({\xi}^{\star})$-valued 1-form $a$ by the formula (2.6).
Then ${\nabla}+{\Gamma}({\alpha})$ is a $spin^c$-connection
on $W$, compatible with ${\nabla}+a$ on ${\xi}^{\star}$.
\begin{flushright}Q.E.D.
\end{flushright}

$\bf {Corollary \: 2.4}$ Suppose ${\nabla}^1$, ${\nabla}^2$
are compatible with the pseudohermitian connection. Then
they differ by an imaginary valued 1-form.

\medskip

Note that in this case, the $a$ in the above proof vanishes.

\bigskip

\section{The Weitzenbock formula and the equations}
\setcounter{equation}{0}

Given a $spin^c$-structure $(W,{\Gamma})$ on the dual
contact bundle ${\xi}^{\star}$ and a $spin^c$-connection
${\nabla}$ on $W$, compatible with the pseudohermitian
connection on ${\xi}^{\star}$. We define the associated
Dirac operator $D_{\xi}$ by

$$
D_{\xi}{\Phi}={\Sigma}_{j=1}^{2}{\Gamma}(e^{j}){\nabla}_{e_{j}}
{\Phi}
$$

\noindent for ${\Phi}$ being a section of $W$ and $\{ e^{j},
j=1,2 \}$ being the dual of an orthonormal basis $\{ e_{j},
j=1,2 \}$ in $\xi$. 

Let $e_0$ or $T$ denote the vector field characterized by 
${\theta}(T)=1$ and ${\cal L}_{T}{\theta}=0$. Define the
divergence $div(v)$ of a vector field $v$ with respect to
the pseudohermitian connection ${\nabla}^{{\psi}.h.}$ by

$$
div(v)=
{\Sigma}_{i=0}^{2}<{\nabla}^{{\psi}.h.}_{e_{i}}v,e^{i}>
$$

\noindent 
(note that $e^{0}={\theta}$, $<,>$ is the pairing, and the
definition is independent of the choice of general bases)
It follows that ${\cal L}_{v}({\theta}{\wedge}d{\theta})
=div(v){\theta}{\wedge}d{\theta}$. So we have

\begin{equation}
{\int}div(v){\theta}{\wedge}d{\theta}=0
\end{equation}

\noindent for $M$ being closed (i.e. compact without
boundary). (hereafter, we'll make this assumption)

Since the $spin^c$-connection $\nabla$ is Hermitian, it
is easy to show by (3.1) that its adjoint ${\nabla}^{\star}$
satisfies the following formula

\begin{equation}
{\nabla}^{\star}_{v}{\Phi}=-{\nabla}_{v}{\Phi}
-div(v){\Phi}
\end{equation}

\noindent for a section $\Phi$ of $W$. Let $D_{\xi}^{\star}$ denote 
the adjoint of $D_{\xi}$. Writing $D_{\xi}^{\star}= 
{\Sigma}_{i=1}^{2}{\nabla}^{\star}_{e_{i}}
({\Gamma}(e^{i}))^{\star}$ and using (2.1) and (3.2), we
obtain $D_{\xi}^{\star}=D_{\xi}$, i.e. $D_{\xi}$ is self-adjoint.
(we may assume ${\nabla}^{{\psi}.h.}_{e_{i}}e^{j}=0$,
hence $div(e_{i})$=0, at a point in the computation [Le1]) 
 
Now we compute

\begin{eqnarray}
D_{\xi}^{\star}D_{\xi}{\Phi}&=&D_{\xi}^{2}{\Phi}\:(D_{\xi}
\:being\:self-adjoint)\nonumber\\
 &=&{\Gamma}(e^{i}){\nabla}_{e_{i}}
({\Gamma}(e^{j}){\nabla}_{e_{j}}{\Phi})\:(summation\:convention)
\nonumber\\
 &=&{\Gamma}(e^{i}){\Gamma}(e^{j}){\nabla}_{e_{i}}{\nabla}_{e_{j}}
{\Phi}\:({\nabla}^{{\psi}.h.}_{e_{i}}e^{j}=0\:at\:a\:point\:p)
\\
 &=&{\nabla}^{\star}_{e_{i}}{\nabla}_{e_{i}}{\Phi}+
{\Sigma}_{i<j}{\Gamma}(e^{i}){\Gamma}(e^{j})
({\nabla}_{e_{i}}{\nabla}_{e_{j}}-{\nabla}_{e_{j}}{\nabla}_{e_{i}})
{\Phi}\nonumber\\
 & & (by\:(3.2)\:evaluated\:at\:p)\nonumber
\end{eqnarray}

It is easy to show from the structural equations of pseudohermitian
geometry that $[e_{1},e_{2}]=-2T$ at $p$. (cf. (5.8) in $\S$5) Using this,
we can rewrite (3.3) as follows:

\begin{eqnarray}
D_{\xi}^{\star}D_{\xi}{\Phi}&=&{\Sigma}_{i=1}^{2}{\nabla}^
{\star}_{e_{i}}{\nabla}_{e_{i}}{\Phi}+ \\
 & & {\Gamma}(e^{1}){\Gamma}(e^{2})F^{\nabla}(e_{1},e_{2}){\Phi}
+{\Gamma}(e^{1}){\Gamma}(e^{2}){\nabla}_{-2T}{\Phi} 
\nonumber
\end{eqnarray}

\noindent where $F^{\nabla}(e_{1},e_{2})=[{\nabla}_{e_{1}},
{\nabla}_{e_{2}}]-{\nabla}_{[e_{1},e_{2}]}$ is the curvature
operator in the directions $e_{1},e_{2}$.

For $(W,{\Gamma})=(W_{can},{\Gamma}_{can})$, we can have more
precise description with respect to $\{ 1,\frac{1}{\sqrt 2}{\theta}^
{\bar 1}\}$, a basis of $W_{can}$. Write $\Phi$ as a colume
vector with respect to this basis:

$$
{\Phi}=\left( \begin{array}{c}
               {\alpha}\\{\beta}_{\bar 1}
              \end{array} \right)
\:for\:{\Phi}={\alpha}+{\beta}_{\bar 1}\frac{1}{\sqrt 2}{\theta}^
{\bar 1}.
$$

By (2.4),(2.5), we can write ${\Gamma}={\Gamma}_{can}$ as
matrices:

$$
{\Gamma}(e^{1})=\left( \begin{array}{cc}
                        0&-1\\1&0 \end{array} \right),\:
{\Gamma}(e^{2})=\left( \begin{array}{cc} 
                        0&i\\i&0 \end{array} \right).
$$

The canonical $spin^c$-connection ${\nabla}_{can}$ has the
connection form: $\left( \begin{array}{cc}
                   0&0\\0&i{\omega} \end{array} \right)$
where $\omega$ is the pseudohermitian connection form:
${\nabla}^{{\psi}.h.}e^{1}={\omega}e^{2}$ as in the proof
of Proposition 2.1. So by Corollary 2.4, our $spin^c$
-connection $\nabla$ (compatible with ${\nabla}^{{\psi}.h.}$)
equals $d+A$ with

$$
A=\left( \begin{array}{cc}
          ia & 0 \\ 0 & i({\omega}+a) \end{array} \right)
$$

\noindent where $a$ is a real-valued 1-form. Let $Z_{1}=
\frac{1}{2}(e_{1}-ie_{2})$. A direct computation shows

\begin{equation}
D_{\xi}{\Phi}=\left( \begin{array}{c}
                -2{\beta}_{{\bar 1},1}^{a}\\
                 2{\alpha}_{,{\bar 1}}^{a} \end{array} \right)
\end{equation}

\noindent in which ${\beta}_{{\bar 1},1}^{a}={\beta}_{{\bar 1},1}
+ia(Z_{1}){\beta}_{\bar 1}$, ${\alpha}_{,{\bar 1}}^{a}= 
{\alpha}_{,{\bar 1}}+ia(Z_{\bar 1}){\alpha}$. (covariant
derivative without upper index ``$a$'' is with respect to
the pseudohermitian connection)
  
Observe that $d{\omega}(e_{1},e_{2})=-2{\cal W}$ where $\cal W$
denotes the Tanaka-Webster curvature. ([CL],[Tan],[Web], or (5.4)
in $\S$5) We compute

\begin{eqnarray}
F^{\nabla}(e_{1},e_{2})&=&dA(e_{1},e_{2})\\
 &=&\left( \begin{array}{cc}
            ida(e_{1},e_{2})&0\\
            0&-2i{\cal W}+ida(e_{1},e_{2}) \end{array} \right)
\nonumber
\end{eqnarray}

Taking the Hermitian inner product with $\Phi$ in (3.4) and using
(3.6), we obtain

\begin{eqnarray}
\|D_{\xi}{\Phi}\|^{2}&=&{\Sigma}_{j=1}^{2}\|{\nabla}_{e_{j}}
{\Phi}\|^{2}+2{\int_M}{\cal W}|{\beta}_{\bar 1}|^{2}dv_{\theta}
+\\
 &{\int_M}&da(e_{1},e_{2})(|{\alpha}|^{2}-
|{\beta}_{\bar 1}|^{2})dv_{\theta}
+2i{\int_M}({\alpha}_{,0}^{a}{\bar {\alpha}}-{\beta}_{{\bar 1},0}
^{a}{\beta}_{1})dv_{\theta}\nonumber 
\end{eqnarray}

\noindent in which $dv_{\theta}={\theta}{\wedge}d{\theta}$.
(here $'',0''$ means the covariant derivative in the
$T$-direction) Define ${\pi}_{\xi}$ from 2-forms to functions
by ${\pi}_{\xi}({\eta})={\eta}(e_{1},e_{2})$, i.e. projecting
$\eta$ onto its $e^{1}{\wedge}e^{2}$-component. It is easy
to see ($tr$ means trace) 

\begin{equation}
\frac{1}{2}{\pi}_{\xi}{\circ}tr(F^{{\nabla}-{\nabla}_{can}})
=ida(e_{1},e_{2})
\end{equation}

Let ${\Phi}^{\sigma}=({\alpha},-{\beta})$ for ${\Phi}=
({\alpha},{\beta}){\in}{\cal C}{\oplus}{\Lambda}^{0,1}
{\xi}^{\star}$. Now we can define our ``monopole'' equations
for $(A,{\Phi})$ as follows:

\begin{equation}
\left\{ \begin{array}{l}
         D_{\xi}{\Phi}=0\\ \frac{1}{2}{\pi}_{\xi}{\circ}
      tr(F^{\nabla})=i<{\Phi}^{\sigma},{\Phi}> 
       _{h} \end{array} \right.
\end{equation}

\noindent in which $<,>_h$ denotes the Hermitian inner product
induced by $h$ on $W_{can}$. Recall that on a $CR$ or 
pseudohermitian manifold, we have ${\bar {\partial}}_b$-operator
mapping ${\Lambda}^{p,q}$ to ${\Lambda}^{p,q+1}$. Also
with respect to the connection ${\nabla}={\nabla}_{can}+ia$,
we have the associated covariant differentiation
${\bar {\partial}}_b^{a}$. For our case, ${\bar {\partial}}_b^{a}
{\alpha}={\alpha}_{,{\bar 1}}^{a}{\theta}^{\bar 1}$ for $\alpha$
being a function while $({\bar {\partial}}_b^{a})^{\star}
{\beta}=-{\sqrt 2}{\beta}_{{\bar 1},1}^{a}$ for ${\beta}=
{\beta}_{\bar 1}\frac{1}{\sqrt 2}{\theta}^{\bar 1}$.
(note that $\frac{1}{\sqrt 2}{\theta}^{\bar 1}$ has length 1
with respect to $<,>_h$)

Now by (3.5) it is clear that 

\begin{equation}
D_{\xi}={\sqrt 2}({\bar {\partial}}_b^{a}+
({\bar {\partial}}_b^{a})^{\star})
\end{equation}

Therefore in terms of $(a,{\alpha},{\beta}={\beta}_{\bar 1}
\frac{1}{\sqrt 2}{\theta}^{\bar 1})$, (3.9)
is equivalent to 

\begin{eqnarray}
\left\{ \begin{array}{c}({\bar {\partial}}_b^{a}+
({\bar {\partial}}_b^{a})^{\star})({\alpha}+{\beta})=0\\
(or\:{\alpha}_{,{\bar 1}}^{a}=0,\:{\beta}_{{\bar 1},1}^{a}=0)\\
da(e_{1},e_{2})-{\cal W}=|{\alpha}|^{2}-|{\beta}_{\bar 1}|^{2}
\end{array} \right.
\end{eqnarray}

\noindent by (3.10), (3.5), and (3.8).
\medskip

$\bf {Proof\:of\:Theorem\:B}$: Substituting (3.11) and (1.1)
in (3.7) gives

\begin{eqnarray}
0 & = & {\Sigma}_{j=1}^{2}\|{\nabla}_{e_{j}}
{\Phi}\|^{2}+{\int_M}{\cal W}(|{\alpha}|^{2}+
|{\beta}_{\bar 1}|^{2})dv_{\theta}
\\
 & & +{\int_M}(|{\alpha}|^{2}-|{\beta}_{\bar 1}|^{2})^{2}
dv_{\theta}\nonumber
\end{eqnarray}

Now the theorem follows from (3.12). 
\begin{flushright}Q.E.D.
\end{flushright}      

\bigskip

\section{Proof of Theorem A} 
\setcounter{equation}{0}




We define an almost complex structure $\tilde J$ on $M{\times}R$, the
``symplectification'' of the contact manifold $(M,{\xi})$ as follows:
${\tilde J}=J$ on $\xi$, ${\tilde J}(e_{3})=e_{0},{\tilde J}(e_{0})
=-e_{3}$. Here $e_{3}={\partial}/{\partial}t$, $t$ being the coordinate
of $R$, and recall that $e_{0}$ is just the vector field $T$. (see
$\S$3 or $\S$5) Let $g=(dt)^{2}+h$ where $h$ is the adapted metric.
($\S$5) Let $\{ e^{j},j=0,1,2,3 \}$ be the dual
basis of the orthonormal basis $\{ e_{j},j=0,1,2,3 \}$ with respect
to the metric $g$. (recall that $e_{j}$ in $\xi$ and $e^{j}$ in
${\xi}^{\star}$ for $j=1,2$ are defined in $\S$2. Of course we have
viewed ${\xi}^{\star}$ as a subset of $T^{\star}(M{\times}R)$)
$\tilde J$ also acts on cotangent vectors by $({\tilde J}v)(w)=
v({\tilde J}w)$ as usual. 
Associated to ${\tilde J}$, we have a canonical $spin^c$-structure on
$(M{\times}R, g)$. The differential forms of type $(0,{\star})$ 
constitute the spinors. The Clifford multiplication is defined by

$$
{\Gamma}(w){\tau}=\frac{1}{\sqrt 2}w''{\wedge}{\tau}
-{\sqrt 2}{\iota}(w_{\#}){\tau}.
$$

\noindent (cf. [Sal], for instance) Here $w_{\#}$ denotes the 
corresponding tangent vector of the cotangent vector $w$ with
respect to $g$, and $w''=w+i{\tilde J}w$. Let ${\theta}^{\bar 2}
=e^{3}-ie^{0}$. It is easy to compute that $(e^{3})''={\theta}^{\bar 2},
(e^{0})''=i{\theta}^{\bar 2}$. (similarly, $(e^{1})''={\theta}^{\bar 1},
(e^{2})''=i{\theta}^{\bar 1}$. Recall that ${\theta}^{\bar 1}=
e^{1}-ie^{2}$) Define a map
 
$$
{\varpi}:{\cal C}{\oplus}
{\Lambda}^{0,2}T^{\star}(M{\times}R){\rightarrow}W_{can}   
={\cal C}{\oplus}{\Lambda}^{0,1}{\xi}^{\star} 
$$

\noindent by deleting the $\frac{1}{\sqrt 2}{\theta}^{\bar 2}$ factor and
restricting its domain of definition. (the first $\cal C$ denotes
the trivial complex line bundle over $M{\times}R$ while the second
$\cal C$ means the trivial complex line bundle over $M$) In practice,
we write ${\varpi}({\alpha}+{\beta_{\bar 1}}{\theta}^{\bar 1}{\wedge}
{\theta}^{\bar 2})={\alpha}+{\sqrt 2}{\beta_{\bar 1}}{\theta}^{\bar 1}$.
Conversely by extending the domain of definition and wedging
$\frac{1}{\sqrt 2}{\theta}^{\bar 2}$ in the second component, we get
a map ${\Xi}:W_{can}{\rightarrow}{\cal C}{\oplus}
{\Lambda}^{0,2}T^{\star}(M{\times}R)$ with ${\varpi}{\circ}{\Xi}$
being the identity. We often write ${\tilde {\Phi}}$ instead of
${\Xi}({\Phi})$.Now we can define ${\rho}:T^{\star}M{\rightarrow}
End(W_{can})$ by

$$
{\rho}(w)({\Phi})={\varpi}{\Gamma}(e^{3}){\Gamma}(w)({\tilde{\Phi}}).
$$

Let $\Phi_0$ be the canonical section $(1,0)$ in $W_{can}$. Let
${\Phi_1}={\frac{1}{\sqrt 2}}{\theta}^{\bar 1}$. A direct computation
shows that ${\rho}(e^{0})({\Phi_0})=-i{\Phi_0}$,
${\rho}(e^{0})({\Phi_1})=i{\Phi_1}$, ${\rho}(e^{1})({\Phi_0})=
-{\Phi_1}$, ${\rho}(e^{1})({\Phi_1})={\Phi_0}$,
${\rho}(e^{2})({\Phi_0})=-i{\Phi_1}$, ${\rho}(e^{2})({\Phi_1})=-i{\Phi_0}$.
In matrix form with respect to the orthonormal basis $\{ {\Phi_0},
{\Phi_1} \}$, we have

$$ 
{\rho}(e^{0})=\left( \begin{array}{cc}
                        -i& 0\\0& i \end{array} \right),\:
{\rho}(e^{1})=\left( \begin{array}{cc} 
                        0&1\\-1&0 \end{array} \right),\:
{\rho}(e^{2})=\left( \begin{array}{cc} 
                        0&-i\\-i&0 \end{array} \right).
$$

Now it is clear that $\rho$ defines a Clifford multiplication.
And from the above construction
\medskip

\noindent $\bullet$ $(W_{can},{\rho})$ is isomorphic to
the restriction to $M$ of the canonical $spin^c$-structure induced by
${\tilde J}$. (see the definition of ``restriction'' in the proof
of Proposition 4.3 in [Kro])
\medskip
 
There is a canonical $spin^c$-connection ${\tilde{\nabla}}^{can}$
on ${\cal C}{\oplus}
{\Lambda}^{0,2}T^{\star}(M{\times}R)$ which is compatible with
the Levi-Civita connection ${\nabla}^g$ of $g$. ([Sal]) We define
a connection ${\nabla}^{can}$ on $W_{can}$ by

$$
{\nabla}^{can}_{v}{\Phi}={\varpi}({\tilde{\nabla}}^{can}_{v}{\tilde{\Phi}})
$$

\noindent for $v$ in $TM{\subset}T(M{\times}R)$. Let ${\nabla}^{h}$
denote the Levi-Civita connection of the metric $h$ on $M$.
Noting that ${\nabla}^{g}
_{v}e^{3}=0,{\nabla}^{g}_{v}w={\nabla}^{h}_{v}w$ for $v$ in $TM$, $w$
in $T^{\star}M$ (viewed as a subset of $T^{\star}(M{\times}R)$),
we can easily verify that ${\nabla}^{can}$ is a $spin^c$-connection
on $(W_{can},{\rho})$, compatible with the Levi-Civita connection
${\nabla}^h$. Note that ${\nabla}^{can}$ is different from 
${\nabla}_{can}$ in $\S$2 which is compatible with the pseudohermitian
connection on ${\xi}^{\star}$. To use the ``monopole class'' condition,
we will choose a special family of Riemannian metrics on $M$. Let

$$
h_{\epsilon}=({\epsilon}e^{0})^{2}+(e^{1})^{2}+(e^{2})^{2}
$$

\noindent and Let $g_{\epsilon}=(dt)^{2}+h_{\epsilon}$ be the
corresponding metric on $M{\times}R$. (recall that $e^{3}=dt$)
So $e^{0}_{\epsilon}={\epsilon}e^{0}$, $e^{1}$, $e^{2}$, (and
$e^{3}$, resp.) form an orthonormal coframe for $h_{\epsilon}$
($g_{\epsilon}$, resp.) Now with $e^{0}_{\epsilon}$, $g_{\epsilon}$,
$h_{\epsilon}$ replacing $e^{0}$, $g$, $h$ resp., we can go through
the above procedure again to
get ${\tilde J}_{\epsilon}$, ${\Gamma}_{\epsilon}$, 
${\theta}^{\bar 2}_{\epsilon}$, $\varpi_{\epsilon}$, $\Xi_{\epsilon}$,
${\tilde{\Phi}}^{\epsilon}$, ${\rho}_{\epsilon}$,
$ ^{\epsilon}{\tilde{\nabla}}^{can}$, and $ ^{\epsilon}{\nabla}^{can}$.
Note that the hermitian metric on $W_{can}$ does not change. It is easy
to verify that ${\rho}_{\epsilon}(e^{j}_{\epsilon})={\rho}(e^{j})$
for $j=0,1,2$. Here $e^{j}_{\epsilon}=e^{0}_{\epsilon}$ if $j=0$;
$=e^{j}$ otherwise. Also $ ^{\epsilon}{\nabla}^{can}$ is a $spin^c$-
connection on $(W_{can}, {\rho_{\epsilon}})$, compatible with the
Levi-Civita connection ${\nabla}^{h_{\epsilon}}$. Recall (see $\S$5,
the Appendix) that ${A^1}_{\bar 1}=
A_{{\bar 1}{\bar 1}}$ ($h_{1{\bar 1}}=1$) denotes the pseudohermitian
torsion with respect to $(J,{\theta})$.
\medskip

$\bf {Proposition \: 4.1}$ (1) $ ^{\epsilon}{\nabla}^{can}{\Phi_0}=
\frac{i}{\sqrt 2}{\epsilon^{-1}}A_{{\bar 1}{\bar 1}}{\theta^1}{\otimes}
{\Phi_1}$.

(2) $ ^{\epsilon}{\nabla}^{can}{\Phi_1}=\frac{i}{4}{\epsilon^{-1}}A_{11}
{\theta^1}{\otimes}{\Phi_0}+i({\omega}+{\epsilon}{\theta})
{\otimes}{\Phi_1}$.
\medskip

$\bf{Proof}$: Let us review how to obtain $ ^{\epsilon}{\tilde{\nabla}}^{can}$
from the Levi-Civita connection $\nabla^{g_{\epsilon}}$ on $M{\times}R$.
([Sal]) Let $\Psi$ be an endomorphism of the tangent bundle. Define
${\iota}({\Psi})$ acting on a $k$-form $\tau$ by

$$
{\iota}({\Psi}){\tau}(v_{1},...,v_{k})={\Sigma}_{j=1}^{k}
{\tau}(v_{1},...,v_{j-1},{\Psi}v_{j},v_{j+1},...,v_{k})
$$

\noindent for tangent vectors $v_{1},...,v_{k}$. (${\iota}({\Psi}){\tau}=0$
if $\tau$ is a function) Let $N_{\epsilon}$ denote the Nijenhuis tensor of
${\tilde J}_{\epsilon}$. Our canonical $spin^c$-connection 
$ ^{\epsilon}{\tilde{\nabla}}^{can}$ is defined by

\begin{equation}
 ^{\epsilon}{\tilde{\nabla}}^{can}_{v}{\tau}={\nabla}^{g_{\epsilon}}_
{v}{\tau}+\frac{1}{2}{\iota}({\tilde J}_{\epsilon}
{\nabla}^{g_{\epsilon}}_{v}{\tilde J}_{\epsilon}){\tau}
+\frac{1}{8}{\Theta_v^{\epsilon}}{\wedge}{\tau}+
\frac{1}{2}{\iota}({\bar{\Theta}}_v^{\epsilon}){\tau}
\end{equation}

\noindent in which ${\Theta_v^{\epsilon}}$ is a (0,2)-form defined by
${\Theta_v^{\epsilon}}(x,y)=g_{\epsilon}(v, N_{\epsilon}(x,y))$, and
${\iota}$ in the last term is just the usual interior product (of forms).

Let $\{ Z_{1},Z_{\bar 1},Z_{2}^{\epsilon},Z_{\bar 2}^{\epsilon}\}$ be
a basis dual to $\{ {\theta^1},{\theta^{\bar 1}},{\theta^{2}_{\epsilon}},
{\theta^{\bar 2}_{\epsilon}} \}$. A direct computation using the
formula $[Z_{\bar 1},T]={A^1}_{\bar 1}Z_{1}-{\omega_{\bar 1}}^{\bar 1}(T)
Z_{\bar 1}$ (cf. (5.9) in $\S$5) shows that
$N_{\epsilon}(Z_{\bar 1},Z_{\bar 2}^{\epsilon})=2i{\epsilon^{-1}}
A_{{\bar 1}{\bar 1}}Z_{1}$.  Hence

\begin{equation}
{\Theta_v^{\epsilon}}=g_{\epsilon}(v,2i{\epsilon^{-1}}
A_{{\bar 1}{\bar 1}}Z_{1}){\theta^{\bar 1}}{\wedge}
{\theta^{\bar 2}_{\epsilon}}
\end{equation}

\noindent (1) follows from (4.1) easily. To compute 
$ ^{\epsilon}{\nabla}^{can}{\Phi_1}$, we need to know $
^{\epsilon}{\tilde{\nabla}}^{can}({\theta^{\bar 1}}{\wedge}
{\theta^{\bar 2}_{\epsilon}})$. Let $\omega_{j({\epsilon})}^{i}$ be
the Riemannian
connection forms for $h_{\epsilon}$ so that ${\nabla^{h_{\epsilon}}}{e^{i}_
{\epsilon}}=-{\omega_{j({\epsilon})}^{i}}{\otimes}{e^{j}_{\epsilon}}$.
Then in the tangent direction of $M$,
we compute

\begin{eqnarray}
{\nabla^{g_{\epsilon}}}{\theta^{\bar 1}}& = &{\nabla^{h_{\epsilon}}}
{\theta^{\bar 1}}=-{\omega_{j({\epsilon})}^{1}}{\otimes}{e^{j}_{\epsilon}}
-i(-{\omega_{j({\epsilon})}^{2}}{\otimes}{e^{j}_{\epsilon}})\\
 &=&i{\omega_{1({\epsilon})}^{2}}{\otimes}(e^{1}-ie^{2})-(
{\omega_{0({\epsilon})}^{1}}-i{\omega_{0({\epsilon})}^{2}}){\otimes}
{e^0_{\epsilon}} \nonumber\\
 &=&i({\omega}+{\epsilon}{\theta}){\otimes}{\theta^{\bar 1}}+
(i{\epsilon}{\theta^{\bar 1}}-{A^{\bar 1}}_{1}{\theta^{1}}){\otimes}
{\theta}\nonumber
\end{eqnarray}

\noindent by (5.6) and (5.7) for the metric $h_{\epsilon}$. Note that 
$e^0_{\epsilon}={\epsilon}{\theta}$ and the torsion $A_{11}^{\epsilon}$
for $(J,e^0_{\epsilon})$ equals ${\epsilon}^{-1}A_{11}$.

Let $\{ e_{j}^{\epsilon},j=0,1,2,3 \}$ denote the basis dual to
$\{ e^{j}_{\epsilon},j=0,1,2,3 \}$. Then it is easy to see that
$e_{0}^{\epsilon}={\epsilon^{-1}}T$, $e_{j}^{\epsilon}=e_{j}$
for $j=1,2,3$, and ${Z_1^{\epsilon}}=Z_{1}=\frac{1}{2}(e_{1}-ie_{2})$, $ 
Z_{2}^{\epsilon}=\frac{1}{2}(e_{3}-ie_{0}^{\epsilon})$.
Using ${\nabla}^{g_{\epsilon}}_{v}e^{3}=0$ for $v$ in $TM$ and (5.6),
(5.7), we can show that

$$
{\nabla^{g_{\epsilon}}}Z_{\bar 1}=-i({\omega}+{\epsilon}{\theta})Z_{\bar 1}
-\frac{1}{2}(i{\theta^1}+{\epsilon^{-1}}{A^1}_{\bar 1}{\theta^{\bar 1}})
{\epsilon^{-1}}T.
$$

\noindent It follows that 

\begin{eqnarray}
({\tilde J}_{\epsilon}{\nabla^{g_{\epsilon}}}{\tilde J}_{\epsilon})Z_{\bar 1} 
& = &(-i){\tilde J}_{\epsilon}{\nabla^{g_{\epsilon}}}Z_{\bar 1}+
{\nabla^{g_{\epsilon}}}Z_{\bar 1}\nonumber\\
 & = & ({\theta^1}-i{\epsilon}^{-1}{A^1}_{\bar 1}{\theta^{\bar 1}}){\otimes}
Z_{2}^{\epsilon}\nonumber
\end{eqnarray}

Similarly using ${\omega}_{0}^{1}+i{\omega}_{0}^{2}=i{\theta^1}+
{A^1}_{\bar 1}{\theta^{\bar 1}}$ (the complex version of (5.7)) for
$h_{\epsilon}$, we can easily obtain

$$
({\tilde J}_{\epsilon}{\nabla^{g_{\epsilon}}}{\tilde J}_{\epsilon})
Z_{\bar 2}^{\epsilon}=(-{\theta^1}+i{\epsilon}^{-1}{A^1}_{\bar 1}
{\theta^{\bar 1}}){\otimes}Z_{1}
$$

Since $[{\iota}({\tilde J}_{\epsilon}{\nabla_{v}^{g_{\epsilon}}}
{\tilde J}_{\epsilon}){\theta^{\bar 1}}](w)={\theta^{\bar 1}}(
({\tilde J}_{\epsilon}{\nabla_{v}^{g_{\epsilon}}}{\tilde J}_{\epsilon})(w))$,
it follows from the above two formulas that

\begin{equation}
{\iota}({\tilde J}_{\epsilon}{\nabla_{v}^{g_{\epsilon}}}{\tilde J}_{\epsilon})
{\theta^{\bar 1}}=(-{\theta^{\bar 1}}-i{\epsilon}^{-1}{A^{\bar 1}}_{1}
{\theta^1})(v){\otimes}{\theta^{2}_{\epsilon}}
\end{equation}

\noindent for $v$ in $TM$. Replacing ${\theta^{\bar 1}}$ by ${\theta}^{\bar 2}
_{\epsilon}$ in the previous computation, we obtain

\begin{eqnarray}
&&{\nabla_{v}^{g_{\epsilon}}}{\theta}^{\bar 2}_{\epsilon}
=-\frac{1}{2}(i{\epsilon}^{-1}{A^{\bar 1}}_{1}{\theta^1}+{\theta^{\bar 1}})
(v){\theta^1}+\frac{1}{2}({\theta^1}-i{\epsilon}^{-1}{A^{1}}_{\bar 1}
{\theta^{\bar 1}})(v){\theta^{\bar 1}} \\
&&{\iota}({\tilde J}_{\epsilon}{\nabla_{v}^{g_{\epsilon}}}{\tilde J}_{\epsilon})
{\theta}^{\bar 2}_{\epsilon}=({\theta^{\bar 1}}+i{\epsilon}^{-1}
{A^{\bar 1}}_{1}{\theta^1})(v){\theta^1}
\end{eqnarray}

\noindent for $v$ in $TM$. On the other hand, it is easy to see that 

\begin{equation}
{\iota}({\bar{\Theta}}^{\epsilon}_{v}){\theta^{\bar 1}}{\wedge}
{\theta}^{\bar 2}_{\epsilon}=i{\epsilon}^{-1}{A^{\bar 1}}_{1}{\theta^1}(v)
\end{equation}

\noindent by (4.2). Let ${\tilde{\nabla}}^{g_{\epsilon}}$ denote
the sum of ${\nabla}^{g_{\epsilon}}$ and $\frac{1}{2}
{\iota}({\tilde J}_{\epsilon}{\nabla^{g_{\epsilon}}}{\tilde J}_{\epsilon})$.
Now we can compute

\begin{eqnarray}
^{\epsilon}{\nabla}^{can}{\Phi_1}&=&\frac{1}{2}{\varpi_{\epsilon}}
^{\epsilon}{\tilde{\nabla}}^{can}({\theta^{\bar 1}}{\wedge}
{\theta}^{\bar 2}_{\epsilon})\nonumber\\
&=&\frac{1}{2}{\varpi_{\epsilon}}[({\tilde{\nabla}}^{g_{\epsilon}}{\theta^{\bar 1}}){\wedge}
{\theta}^{\bar 2}_{\epsilon}+{\theta^{\bar 1}}{\wedge}(
{\tilde{\nabla}}^{g_{\epsilon}}{\theta}^{\bar 2}_{\epsilon})]+
\frac{1}{4}{\iota}({\bar{\Theta}}^{\epsilon}_{}){\theta^{\bar 1}}{\wedge}
{\theta}^{\bar 2}_{\epsilon}\nonumber\\
&=&i({\omega}+{\epsilon}{\theta}){\otimes}{\Phi_1}+\frac{1}{4}
i{\epsilon}^{-1}{A^{\bar 1}}_{1}{\theta^1}{\otimes}{\Phi_0}
\nonumber
\end{eqnarray}

\noindent by (4.3),(4.4),(4.5),(4.6), and (4.7).

\begin{flushright}Q.E.D.
\end{flushright}

Next we'll deal with the Dirac operator $D_{A_{\epsilon}}$ associated to 
the canonical $spin^c$-connection $^{\epsilon}{\nabla}^{can}$.
Here $A_{\epsilon}$ denotes the connection form with respect to the basis
$\{ {\Phi_0},{\Phi_1} \}$:

$$
\left( \begin{array}{cc}
       0 & \frac{i}{4}{\epsilon^{-1}}A_{11}{\theta^1}\\
      \frac{i}{\sqrt 2}{\epsilon^{-1}}A_{{\bar 1}{\bar 1}}{\theta^1}
       & i({\omega}+{\epsilon}{\theta})
       \end{array} \right).
$$

The Clifford multiplication $\rho_{\epsilon}$ of ${\eta}=de^{0}=2e^{1}{\wedge}e^{2}$
can be easily computed:

\begin{eqnarray} 
{\rho_{\epsilon}}({\eta}){\Phi_0}& = &2{\rho_{\epsilon}}
(e^{1}){\rho_{\epsilon}}(e^{2}){\Phi_0}=-2i{\Phi_0}\\
{\rho_{\epsilon}}({\eta}){\Phi_1}& = &2{\rho_{\epsilon}}(e^{1})
{\rho_{\epsilon}}(e^{2}){\Phi_1}=2i{\Phi_1}\nonumber
\end{eqnarray} 

Let $\star_{\epsilon}$ denote the Hodge star-operator with respect to
the metric $h_{\epsilon}$. Since ${\rho_{\epsilon}}(e^{j}_{\epsilon})
{\rho_{\epsilon}}({\Omega})={\rho_{\epsilon}}(e^{j}_{\epsilon}{\wedge}
{\Omega}-{\iota}(e_{j}^{\epsilon}){\Omega})$ for an arbitrary function
or form $\Omega$,
we can compute that for a scalar function or forms $\gamma$

\begin{eqnarray}
{\rho_{\epsilon}}(e^{j}_{\epsilon}){\rho_{\epsilon}}(
\nabla_{e_{j}^{\epsilon}}^{h_{\epsilon}}{\gamma})&=&
{\rho_{\epsilon}}(e^{j}_{\epsilon}{\wedge}\nabla_{e_{j}^{\epsilon}}^{h_{\epsilon}}{\gamma}-{\iota}(e_{j}^{\epsilon})\nabla_{e_{j}^{\epsilon}}^{h_
{\epsilon}}{\gamma})\\
 &=& {\rho_{\epsilon}}((d+d^{\star_{\epsilon}}){\gamma}).
\nonumber
\end{eqnarray}

Note that $d^{\star_{\epsilon}}={\star_{\epsilon}}d{\star_{\epsilon}}$ on
2-forms (changes sign on 1-forms). So for ${\eta}=de^{0}=
2e^{1}{\wedge}e^{2}$, we have

\begin{eqnarray}
d^{\star_{\epsilon}}{\eta} & = & 2{\star_{\epsilon}}(de^{0}_{\epsilon})\\
 & = & 4{\epsilon}{\star_{\epsilon}}(e^{1}{\wedge}e^{2})=4
 {\epsilon}^{2}e^{0}.\nonumber
\end{eqnarray}

Now we can compute $D_{A_{\epsilon}}{\Phi_0}$ as follows:

\begin{eqnarray}
-2iD_{A_{\epsilon}}{\Phi_0}&=&D_{A_{\epsilon}}({\rho_{\epsilon}}({\eta})
{\Phi_0})\:(by\:(4.8))\nonumber\\
& = & {\Sigma}_{j=0}^{2}{\rho_{\epsilon}}(e^{j}_{\epsilon})[{\rho_{\epsilon}}
({\nabla}_{e_{j}^{\epsilon}}^{h_{\epsilon}}{\eta}){\Phi_0}+
{\rho_{\epsilon}}({\eta}){ ^{\epsilon}{\nabla}^{can}_{e_{j}^{\epsilon}}}
{\Phi_0}]\nonumber\\
& = & {\rho_{\epsilon}}((d+d^{\star_{\epsilon}}){\eta}){\Phi_0}
+2iD_{A_{\epsilon}}{\Phi_0}\:(by\:(4.9),Prop.4.1(1),and\:(4.8))
\nonumber\\
& = & -4i{\epsilon}{\Phi_0}+2iD_{A_{\epsilon}}{\Phi_0}\:
(by\:(4.10)\:and\:d{\eta}=0).\nonumber
\end{eqnarray}

\noindent Therefore we obtain

\begin{equation}
D_{A_{\epsilon}}{\Phi_0}={\epsilon}{\Phi_0}.
\end{equation}

Before computing $D_{A_{\epsilon}}{\Phi}$ for a general section $\Phi$
we need two more preparatory formulas. Let $\alpha$ be a scalar function.
It follows easily from (4.9) that

\begin{equation}
{\Sigma}_{j=0}^{2}{\rho_{\epsilon}}(e^{j}_{\epsilon})
{ ^{\epsilon}}{\nabla}_{e_{j}^{\epsilon}}^{can}({\rho_{\epsilon}}
({\alpha}){\Phi_0})={\rho_{\epsilon}}(d{\alpha}){\Phi_0}+{\alpha}
D_{A_{\epsilon}}{\Phi_0}.
\end{equation}

\noindent Also a direct computation shows

\begin{eqnarray}
{\rho_{\epsilon}}({\theta^1}){\Phi_0}=0,\:{\rho_{\epsilon}}
({\theta^{\bar 1}}){\Phi_0}=-2{\Phi_1}\\
{\rho_{\epsilon}}({\theta}{\wedge}{\theta^1}){\Phi_0}=0\nonumber\\
{\rho_{\epsilon}}({\theta}{\wedge}{\theta^{\bar 1}}){\Phi_0}=
-2i{\epsilon}^{-1}{\Phi_1}\nonumber
\end{eqnarray}

Let ${\Phi}={\alpha}{\Phi_0}+{\beta}_{\bar 1}{\Phi_1}$ be a section
of $W_{can}$. (recall ${\Phi_1}=\frac{\theta^{\bar 1}}{\sqrt 2}$)
Under the condition $A_{11}=0$, $^{\epsilon}{\nabla}^{can}{\Phi_0}=0$
by Proposition 4.1(1). We compute, under this condition,

\begin{eqnarray}
D_{A_{\epsilon}}{\Phi}&=&D_{A_{\epsilon}}[({\rho_{\epsilon}}({\alpha})
-\frac{1}{2}{\beta}_{\bar 1}{\rho_{\epsilon}}({\theta^{\bar 1}})){\Phi_0}]
\:(by\:(4.13))\nonumber\\
&=&{\rho_{\epsilon}}(d{\alpha}){\Phi_0}+{\alpha}D_{A_{\epsilon}}{\Phi_0}
-\frac{1}{2}{\rho_{\epsilon}}((d+d^{\star_{\epsilon}})({\beta}_{\bar 1}
{\theta^{\bar 1}})){\Phi_0}\:(by\:(4.12),(4.9))\nonumber\\
&=&-i{\epsilon}^{-1}{\alpha}_{,0}{\Phi_0}-2{\alpha}_{,{\bar 1}}{\Phi_1}
+{\alpha}{\epsilon}{\Phi_0}+{\beta}_{{\bar 1},1}{\Phi_0}+
i{\epsilon}^{-1}{\beta}_{{\bar 1},0}{\Phi_1}+{\beta}_{{\bar 1},1}{\Phi_0}
\nonumber\\
& &(by\:(4.11),(4.13)\:and\:{\star_{\epsilon}}{\theta^{\bar 1}}=i
{\theta^{\bar 1}}{\wedge}e^{0}_{\epsilon})\nonumber\\
&=&(2{\beta}_{{\bar 1},1}-i{\epsilon}^{-1}{\alpha}_{,0}+{\epsilon}{\alpha})
{\Phi_0}+(i{\epsilon}^{-1}{\beta}_{{\bar 1},0}-2{\alpha}_{,{\bar 1}}){\Phi_1}
\nonumber
\end{eqnarray}

It is known that any two $spin^c$-connections compatible with the Levi-Civita
connection differ by an imaginary valued 1-form. (e.g. [Sal])
So we can assume a 
general $spin^c$-connection (on $W_{can}$) compatible with $\nabla^{h_
{\epsilon}}$ has the connection form $A_{\epsilon}+iaI$ (with respect to
the basis $\{ {\Phi_0},{\Phi_1} \}$) with $a$ being a real valued 1-form
and $I$ being a $2{\times}2$ identity matrix. Now we compute

\begin{eqnarray}
D_{A_{\epsilon}+iaI}{\Phi}&=&D_{A_{\epsilon}}{\Phi}+{\Sigma}_{j=0}^{2}
{\rho_{\epsilon}}(e^{j}_{\epsilon})(ia(e_{j}^{\epsilon}){\Phi})
\nonumber\\
&=&(2{\beta}^{a}_{{\bar 1},1}-i{\epsilon}^{-1}{\alpha}^{a}_{,0}+
{\epsilon}{\alpha})
{\Phi_0}+(i{\epsilon}^{-1}{\beta}^{a}_{{\bar 1},0}-2{\alpha}^{a}_
{,{\bar 1}}){\Phi_1}\nonumber
\end{eqnarray}

\noindent in which ${\alpha}^{a}_{,0}={\alpha}_{,0}+ia(T){\alpha}$,
${\alpha}^{a}_{,{\bar 1}}={\alpha}_{,{\bar 1}}+ia(Z_{\bar 1}){\alpha}$,
${\beta}^{a}_{{\bar 1},0}={\beta}_{{\bar 1},0}+ia(T){\beta}_{\bar 1}$,
${\beta}^{a}_{{\bar 1},1}={\beta}_{{\bar 1},1}+ia(Z_{1}){\beta}_{\bar 1}$.
So the Dirac equation $D_{A_{\epsilon}+iaI}{\Phi}=0$ is equivalent to

\begin{equation}
\{ \begin{array}{l}
    2{\beta}^{a}_{{\bar 1},1}-i{\epsilon}^{-1}{\alpha}^{a}_{,0}+
{\epsilon}{\alpha}=0  \\
    i{\epsilon}^{-1}{\beta}^{a}_{{\bar 1},0}-2{\alpha}^{a}_
{,{\bar 1}}=0 
   \end{array}
\end{equation}

Next we'll express the second one of Seiberg-Witten monopole equations
in a workable form. Let $b=\frac{1}{2}tr(A_{\epsilon}+iaI)$. It follows
from Proposition 4.1 that

\begin{equation}
b=\frac{1}{2}i({\omega}+{\epsilon}{\theta})+ia.
\end{equation}

Let $F_A$ denote the curvature 2-form of $A$. Write $F_{b}=iF_{b}^{12}
e^{1}{\wedge}e^{2}+iF_{b}^{01}e^{0}{\wedge}e^{1}+iF_{b}^{02}e^{0}{\wedge}e^{2}
$. It is easy to see

\begin{equation}
{\rho}_{\epsilon}(F_{b})=\left( \begin{array}{cc}
        F_{b}^{12} & {\epsilon}^{-1}(F_{b}^{01}-iF_{b}^{02})\\
        {\epsilon}^{-1}(F_{b}^{01}+iF_{b}^{02}) & -F_{b}^{12}
       \end{array} \right)
\end{equation}

\noindent with respect to the orthonormal basis $\{ {\Phi_0},{\Phi_1} \}$.
On the other hand, the trace free part of the endomorphism ${\Phi}{\otimes}
{\Phi^{\star}}=h({\Phi},{\cdot}){\Phi}$, denoted $\{ {\Phi}{\otimes} 
{\Phi^{\star}} \}$, reads

\begin{equation}
\{ {\Phi}{\otimes}{\Phi^{\star}} \}=\left( \begin{array}{cc}
        \frac{1}{2}(|{\alpha}|^{2}-|{\beta}_{\bar 1}|^{2}) & {\alpha}
        {\beta_{1}}\\
        {\bar{\alpha}}{\beta}_{\bar 1} & \frac{1}{2}(|{\beta}_{\bar 1}|^{2}
        -|{\alpha}|^{2})     
      \end{array} \right)  
\end{equation}

\noindent with respect to the orthonormal basis $\{ {\Phi_0},{\Phi_1} \}$.
(${\beta_1}={\bar{({\beta}_{\bar 1})}}$) From (4.16), (4.17)
the equation ${\rho}_{\epsilon}(F_{b})={\rho}_{\epsilon}
(\frac{1}{2}trF_{A_{\epsilon}+iaI})=\{ {\Phi}{\otimes}{\Phi^{\star}} \}$
is equivalent to the following system:

\begin{equation}
\left\{ \begin{array}{l}
        F_{b}^{12}=\frac{1}{2}(|{\alpha}|^{2}-|{\beta}_{\bar 1}|^{2}) \\
      {\epsilon}^{-1}(F_{b}^{01}+iF_{b}^{02})={\bar{\alpha}}{\beta}_{\bar 1} 
     \end{array} \right.
\end{equation}

Before analyzing the behavior of solutions for the Seiberg-Witten monopole
equations (4.14), (4.18) as ${\epsilon}{\rightarrow}0$, we
need one more result which relates the scalar curvature $R_{h_{\epsilon}}$
of the metric $h_{\epsilon}$ to the Tanaka-Webster curvature $\cal W$ of
the background pseudohermitian structure $(J,{\theta})$. 
\medskip

$\bf{Lemma\:4.2}$: $R_{h_{\epsilon}}=4{\cal W}-{\epsilon}^{2}
-{\epsilon}^{-2}|A_{11}|^{2}$.
\medskip

$\bf{Proof}$: We use the notation in [CH]. Consider a new coframe
${\tilde{\omega}}_{3}={\epsilon}^{2}{\omega_3},{\tilde{\omega}}_{1}=
{\epsilon}{\omega}_{1},{\tilde{\omega}}_{2}={\epsilon}{\omega}_{2}$.
The corresponding connection forms in the structural equations for
the adapted metric ${\epsilon}^{2}h_{\epsilon}=({\tilde{\omega}}_{3})^{2}
+({\tilde{\omega}}_{1})^{2}+({\tilde{\omega}}_{2})^{2}$ read
${\tilde{\psi}}_{3}={\psi}_{3},{\tilde{\psi}}_{1}={\epsilon}^{-1}
{\psi}_{1},{\tilde{\psi}}_{2}={\epsilon}^{-1}{\psi}_{2}$. (note that
${\omega_3},{\omega}_{1},{\omega}_{2}$ are just $e^{0},e^{1},e^{2}$ in
our paper, respectively)

To satisfy (36) in [CH], the $L_{ij}$'s transform as below:
${\tilde L}_{ij}={\epsilon}^{-4}L_{ij}$ for $i,j$ in $\{ 1,2 \}$;
${\tilde L}_{ij}={\epsilon}^{-3}L_{ij}$ if one of indices is $3$.
To determine ${\tilde L}_{33}$ we group the coefficients of
${\omega}_{1}{\wedge}{\omega}_{2}$ in the right-hand side of the third
equation in (36) of [CH] to get

\begin{equation}
\frac{1}{2}{\epsilon}^{-2}(L_{11}+L_{22})+{\epsilon}^{2}{\tilde L}_{33}
=4{\cal W}.
\end{equation}

\noindent (here we have used $d{\psi}_{3}=4{\cal W}{\omega_1}{\wedge}
{\omega_2}$ and note that ${\psi}_{3}$ is just $-{\omega}$ in our
notation) Now we can compute the scalar curvature of the metric
${\epsilon}^{2}h_{\epsilon}$:

\begin{eqnarray}
R_{{\epsilon}^{2}h_{\epsilon}}&=&{\tilde L}_{11}+{\tilde L}_{22}+
{\tilde L}_{33}-1\:([CH])\nonumber\\
 &=&{\epsilon}^{-2}(4{\cal W}-{\epsilon}^{2}-{\epsilon}^{-2}|A_{11}|^{2}).
\nonumber\\
 & &(by\:(4.19)\:and\:\frac{L_{11}+L_{22}}{2}=-|A_{11}|^{2}\:due\:to\:
(38),(40)\:in\:[CH])\nonumber
\end{eqnarray}

Our result follows from the above formula and the dilation relation:
$R_{h_{\epsilon}}={\epsilon}^{2}R_{{\epsilon}^{2}h_{\epsilon}}$.
\begin{flushright} Q.E.D. \end{flushright}

\noindent $\bf{PROOF\:OF\:THEOREM\:A}$: 
\medskip

According to Corollary 5.7 and the proof of Proposition 4.3 in [Kro],
the contact structure $\xi$ being symplectically semifillable implies
that the Euler class $e({\xi})$ is a monopole class for the restriction
to $M$
of the canonical $spin^c$-structure of ``bounded'' symplectic 4-manifold.
The $(W_{can},{\rho_{\epsilon}})$ provides such
a $spin^c$-structure. (note that they are isomorphic to each other for
different $\epsilon$'s and the first Chern class of $W_{can}$ is just
$e({\xi})$) So for the given metric $h_{\epsilon}$, we have a solution
$({\Phi}={\Phi_{\epsilon}},a=a_{\epsilon})$ of (4.14) and (4.18). Recall
that ${\Phi}={\alpha}{\Phi_0}+{\beta}_{\bar 1}{\Phi_1}$, and we sometimes
write ${\alpha}_{\epsilon}$,${\beta}_{\bar 1}^{\epsilon}$ instead of
$\alpha$,${\beta}_{\bar 1}$ to indicate the $\epsilon$-dependence.
 
Now an application of the Weitzenbock formula for the Seiberg-Witten
equations ([Kro] or [Sal]) gives the following estimate:
${\Phi}{\equiv}0$ or, under the assumption $A_{11}=0$,

\begin{equation}
sup|{\Phi}|_{h}^{2}{\leq}sup(-R_{h_{\epsilon}})=sup(-4{\cal W})+
{\epsilon}^{2}
\end{equation}

\noindent by Lemma 4.2. The situation ${\Phi}{\equiv}0$ is ruled out
by the assumption that $e({\xi})$ is not a torsion class: 
${\Phi}{\equiv}0$ implies $F_{b}=0$ which represents the first Chern
class $c_{1}(W_{can})$ of $W_{can}$ up to a constant. But $c_{1}(W_{can})$
is just $e({\xi})$.

The (4.20) tells that $\alpha$ and $\beta_{\bar 1}$ are uniformly
bounded (i.e. there is an upper bound independent of $\epsilon$).
We'll use $O(1)$ to mean an uniformly bounded function or form. Also we
use $O({\epsilon}^{k})$ to mean a function or form bounded by a constant
times ${\epsilon}^{k}$. By (4.18) we have

\begin{equation}
F_{b}^{12}=O(1),F_{b}^{01}=O({\epsilon}),F_{b}^{02}=O({\epsilon}).
\end{equation}

From (4.21) and a theorem of Uhlenbeck (e.g.[Sal]), $b$ is uniformly 
bounded in the $L_{1}^{p}$-norm for any $p>1$ in Coulomb gauges.
(all our norms and the star-operator are with respect to the fixed
metric $h$)
It follows from (4.19) that 
\medskip

$\bf{Lemma\:4.3}$: For a sequence ${\epsilon_j}{\rightarrow}0$,
$a_{\epsilon_j}$ converges weakly in $L_{1}^{p}{\subset}C^{\alpha}$
to $\hat a$.
\medskip

On the other hand we write (4.14) in a matrix form as follows:

\begin{equation}
({\epsilon}^{-1}{\nabla}_{\cal T}^{\epsilon}+{\nabla}_{\Xi}^{\epsilon})
{\Phi}=-{\Lambda}_{\epsilon}{\Phi}
\end{equation}

\noindent in which

$$
{\nabla}_{\cal T}^{\epsilon}=\left( \begin{array}{cc}
            -i{\nabla}_{T}^{a} & 0 \\
            0 & i{\nabla}_{T}^{a} 
           \end{array} \right) ,
{\nabla}_{\Xi}^{\epsilon}=\left( \begin{array}{cc}
            0 & 2{\nabla}_{Z_1}^{a} \\
            -2{\nabla}_{Z_{\bar 1}}^{a} & 0
           \end{array} \right) ,
{\Lambda}_{\epsilon}=\left( \begin{array}{cc}
            {\epsilon} & 0 \\
            0 & 0
           \end{array} \right) .
$$

Taking the square $L^2$-norm of both sides of (4.22) and noting
that ${\nabla}_{\cal T}^{\epsilon},{\nabla}_{\Xi}^{\epsilon},
{\Lambda}_{\epsilon}$ are all self-adjoint, we obtain

\begin{eqnarray}
&&{\epsilon}^{-2}\|{\nabla}_{\cal T}^{\epsilon}{\Phi}\|_{L^2}^{2}
+\|{\nabla}_{\Xi}^{\epsilon}{\Phi}\|_{L^2}^{2}+{\epsilon}^{-1}
<\{{\nabla}_{\cal T}^{\epsilon},{\nabla}_{\Xi}^{\epsilon}\}{\Phi},{\Phi}>
\\ && =<{\Lambda}_{\epsilon}^{2}{\Phi},{\Phi}>={\epsilon}^{2}\|{\alpha}
\|_{L^2}^{2}.\nonumber
\end{eqnarray}

\noindent where $<{\cdot},{\cdot}>$ denotes the $L^2$-inner product 
induced by the metric $h$ and $\{{\nabla}_{\cal T}^{\epsilon},
{\nabla}_{\Xi}^{\epsilon}\}={\nabla}_{\cal T}^{\epsilon}{\nabla}_{\Xi}
^{\epsilon}+{\nabla}_{\Xi}^{\epsilon}{\nabla}_{\cal T}^{\epsilon}$.
\medskip

$\bf{Lemma\:4.4}$: Let $F_{b}^{0}=F_{b}^{01}+iF_{b}^{02}$. Then

\begin{eqnarray}
\{ {\nabla}_{\cal T}^{\epsilon},{\nabla}_{\Xi}^{\epsilon} \} =
\left( \begin{array}{cc}
       0 & 2iA_{11}{\nabla}_{Z_{\bar 1}} \\
       2iA_{{\bar 1}{\bar 1}}{\nabla}_{Z_1} & 0 
       \end{array} \right) +     & & \nonumber \\
\left( \begin{array}{cc}
       0 & {\bar F}_{b}^{0}+iA_{11,{\bar 1}}-2A_{11}a(Z_{\bar 1})\\
       F_{b}^{0}+iA_{{\bar 1}{\bar 1},1}-2A_{{\bar 1}{\bar 1}}a(Z_{1})
       & 0 
     \end{array} \right) . & & \nonumber
\end{eqnarray}

$\bf{Proof}$: A direct computation shows

\begin{equation}
\{{\nabla}_{\cal T}^{\epsilon},{\nabla}_{\Xi}^{\epsilon}\}
\left( \begin{array}{c}
{\alpha}\\{\beta}_{\bar 1} \end{array} \right)
=i\left(\begin{array}{c} 2({\beta}_{{\bar 1},01}^{a}-
{\beta}_{{\bar 1},10}^{a}) \\ 2({\alpha}_{,0{\bar 1}}^{a}-
{\alpha}_{,{\bar 1}0}^{a})
        \end{array} \right) .
\end{equation}

Using the commutation relations: ${\alpha}_{,0{\bar 1}}-{\alpha}_{,{\bar 1}0}
=A_{{\bar 1}{\bar 1}}{\alpha}_{,1}$ and ${\beta}_{{\bar 1},01}-
{\beta}_{{\bar 1},10}={\beta}_{{\bar 1},{\bar 1}}A_{11}+{\beta}_{\bar 1}
A_{11,{\bar 1}}$ ([Le2]), we can compute

\begin{eqnarray}
&&{\beta}_{{\bar 1},01}^{a}-{\beta}_{{\bar 1},10}^{a}
={\beta}_{{\bar 1},{\bar 1}}A_{11}+{\beta}_{\bar 1}A_{11,{\bar 1}}+
i(a_{0,1}-a_{1,0}){\beta}_{\bar 1}\\
&&{\alpha}_{,0{\bar 1}}^{a}-{\alpha}_{,{\bar 1}0}^{a}=
A_{{\bar 1}{\bar 1}}{\alpha}_{,1}+i(a_{0,{\bar 1}}-a_{{\bar 1},0}){\alpha}
\nonumber
\end{eqnarray}

\noindent Here $a_{0}=a(T),a_{1}=a(Z_{1}),a_{\bar 1}=a(Z_{\bar 1})$.
By (4.15) and (5.3) we can easily obtain

\begin{equation}
a_{{\bar 1},0}-a_{0,{\bar 1}}=\frac{1}{2}(F_{b}^{0}+iA_{{\bar 1}{\bar 1},1})
-a_{1}A_{{\bar 1}{\bar 1}}.
\end{equation}

Now Lemma 4.4 follows from (4.24),(4.25),(4.26).
\begin{flushright} Q.E.D.
\end{flushright}

Applying our assumption $A_{11}=0$ and (4.18) to Lemma 4.4 and substituting
the result in (4.23), we obtain

\begin{equation}
{\epsilon}^{2}\|{\alpha}\|_{L^2}^{2}={\epsilon}^{-2}
\|{\nabla}_{\cal T}^{\epsilon}{\Phi}\|_{L^2}^{2}
+\|{\nabla}_{\Xi}^{\epsilon}{\Phi}\|_{L^2}^{2}+ 
2\|{\alpha}{\beta}_{1}\|_{L^2}^{2} 
\end{equation}

\noindent where ${\beta}_{1}$
is the complex conjugate of ${\beta}_{\bar 1}$. It follows that

$$
\|{\alpha}{\beta}_{1}\|_{L^2}^{2}=O({\epsilon}^{2}).
$$

Substituting this in (4.27), we obtain 

\begin{equation}
\|{\nabla}_{\cal T}^{\epsilon}{\Phi}\|_{L^2}^{2}=O({\epsilon}^{4}),
\|{\nabla}_{\Xi}^{\epsilon}{\Phi}\|_{L^2}^{2}=O({\epsilon}^{2}).
\end{equation}

Let ${\hat{\nabla}}_{\cal T}$, ${\hat{\nabla}}_{\Xi}$ denote the
following operators:

\begin{eqnarray}
{\hat{\nabla}}_{\cal T}=\left( \begin{array}{cc}
           -i{\nabla}_{T}^{\hat a} & 0 \\
           0 & i{\nabla}_{T}^{\hat a} 
         \end{array} \right) ,
{\hat{\nabla}}_{\Xi}=\left( \begin{array}{cc}
           0 & 2{\nabla}_{Z_1}^{\hat a} \\
           -2{\nabla}_{Z_{\bar 1}}^{\hat a} & 0
         \end{array} \right) . \nonumber
\end{eqnarray}

It is easy to see that ${\hat{\nabla}}={\hat{\nabla}}_{\cal T}+
{\hat{\nabla}}_{\Xi}$ is an elliptic operator. (independent of 
$\epsilon$) So we can compute

\begin{eqnarray}
\|{\Phi}\|_{L_{1}^2}&{\leq}&C_{1}(\|{\hat{\nabla}}{\Phi}\|_{L^2}+
                            \|{\Phi}\|_{L^2}\:(elliptic\:estimate) \\
&{\leq}&C_{1}(\|{\nabla}_{\cal T}^{\epsilon}{\Phi}\|_{L^2}+
\|{\nabla}_{\Xi}^{\epsilon}{\Phi}\|_{L^2}+\|({\hat{\nabla}}_{\cal T}-
{\nabla}_{\cal T}^{\epsilon}){\Phi}\|_{L^2} \nonumber \\
& & +\|({\hat{\nabla}}_{\Xi}-{\nabla}_{\Xi}^{\epsilon}){\Phi}\|_{L^2}
    +\|{\Phi}\|_{L^2}) \nonumber \\
&{\leq}&C_{2}\:(by\:(4.28),(4.20)) \nonumber
\end{eqnarray}

\noindent in which $C_{1},C_{2}$ are constants independent of 
$\epsilon$, and we can use the covariant derivative ${\nabla}^h$ to
define the Sobolev norm $L_{1}^2$. By (4.29) ${\Phi}={\Phi}_{\epsilon}$
(indicating the $\epsilon$-dependence) converges strongly in $L^2$ for
some sequence ${\epsilon}_j$ tending to $0$. Moreover, 
applying the first inequality of (4.29) to
${\Phi}_{\epsilon_j}-{\Phi}_{\epsilon_k}$ and using (4.28) for 
${\epsilon_j},{\epsilon_k}$ to show $\|{\hat{\nabla}}{\Phi}_{\epsilon_j}
\|_{L^2}$ and $\|{\hat{\nabla}}{\Phi}_{\epsilon_k}\|_{L^2}$ are small
as ${\epsilon_j},{\epsilon_k}$ are small enough, we conclude that
${\Phi}_{\epsilon_j}$ is Cauchy in $L_{1}^{2}$. Therefore
${\Phi}_{\epsilon_j}$ converges strongly in $L_{1}^{2}$ to ${\hat{\Phi}}$
($\alpha_{\epsilon}$, ${\beta}_{\bar 1}^{\epsilon}$ converge to 
${\hat{\alpha}}$, ${\hat{\beta}}_{\bar 1}$, resp.) as ${\epsilon_j}$
goes to $0$. It follows from (4.28) that

\begin{eqnarray}
&{\hat{\nabla}}_{\cal T}{\hat{\Phi}}=0,{\hat{\nabla}}_{\Xi}{\hat{\Phi}}=0&\\
&i.e.\:{\hat{\alpha}}_{,{\bar 1}}^{\hat a}={\hat{\beta}}_{{\bar 1},1}^
{\hat a}=0,\:{\hat{\alpha}}_{,0}^{\hat a}={\hat{\beta}}_{{\bar 1},0}^
{\hat a}=0.& \nonumber
\end{eqnarray}

We'll show the $C^{\infty}$-smoothness of $\hat a$ and $\hat{\Phi}$ by
the usual bootstrap argument. First ${\hat a}{\in}L_{1}^{p}$ ($p>1$) and
${\hat{\Phi}}{\in}L_{1}^{2}$ imply ${\hat a}{\hat{\Phi}}{\in}L_{1}^{2}$
since $L_{1}^{2}{\times}L_{1}^{4}{\subset}L_{1}^{2}$ in dimension 3.
It follows that ${\hat{\Phi}}{\in}L_{2}^{2}$ by the elliptic
regularity. (${\hat{\nabla}}{\hat{\Phi}}=0$ by (4.30))
Since $L_{k}^{2}$ is an algebra for $2k>dimension=3$, $F_{\hat a}$ is
in $L_{2}^{2}$ by (4.18), (4.15). So ${\hat a}$ is in $L_{3}^{2}$. (note
that $d^{\star}{\hat a}=-\frac{1}{2}d^{\star}{\omega}$ is smooth by
(4.15) and $b_{\epsilon}$ having been taken in Coulomb gauges)
Now repeating the above argument, we obtain ${\hat a}{\hat{\Phi}}{\in}
L_{2}^{2}$, then ${\hat{\Phi}}{\in}L_{3}^{2}$,$F_{\hat a}{\in}L_{3}^{2}$,
and ${\hat a}{\in}L_{4}^{2}$, etc.. So ${\hat a}$, ${\hat{\Phi}}$ are
$C^{\infty}$ smooth.
 
On the other hand, taking the limit of the first equation of (4.18) gives

\begin{equation}
d{\hat a}(e_{1},e_{2})-{\cal W}=\frac{1}{2}(|{\hat{\alpha}}|^{2}-
|{\hat{\beta}}_{\bar 1}|^{2}).
\end{equation}

From (4.30),(4.31), $(\frac{\hat{\alpha}}{\sqrt 2},\frac
{{\hat{\beta}}_{\bar 1}}{\sqrt 2},{\hat a})$ is a ($C^{\infty}$ smooth)
solution of (3.11). Suppose both ${\hat{\alpha}}$ and ${\hat{\beta}}_{\bar 1}$
are identically zero. Then by (4.18), $c_{1}(W_{can})=e({\xi})$ vanishes in
$H^{2}(M,C)$, contradicting our assumption.
\bigskip

\section{Appendix: a brief introduction to pseudohermitian geometry}
\setcounter{equation}{0}

Let $M$ be a smooth (paracompact) 3-manifold with an oriented
contact structure
$\xi$. We say a contact form $\theta$ is oriented if $d{\theta}(u,v)>0$
for $(u,v)$ being an oriented basis of $\xi$.
There always exists a global oriented
contact form $\theta$, obtained by patching together local ones with
a partition of unity. The characteristic vector field of $\theta$ is the 
unique vector field $T$ such that ${\theta}(T)=1$ and ${\cal L}_{T}{\theta}
=0$ or $d{\theta}(T,{\cdot})=0$. A $CR$-structure compatible with $\xi$
is a smooth endomorphism $J:{\xi}{\rightarrow}{\xi}$ such that $J^{2}=-
identity$. We say $J$ is oriented if $(X,JX)$ is an oriented basis of
$\xi$ for any nonzero $X{\in}{\xi}$.
A pseudohermitian structure compatible with $\xi$ is a
$CR$-structure $J$ compatible with $\xi$ together
with a global contact form $\theta$.   

Given a pseudohermitian structure $(J,{\theta})$, 
we can choose a complex vector field $Z_1$, an eigenvector
of $J$ with eigenvalue $i$, and a complex 1-form ${\theta}^1$ such that
$\{ {\theta},{\theta^1},{\theta^{\bar 1}} \}$ is dual to $\{ T,Z_{1},Z_
{\bar 1} \}$. (${\theta^{\bar 1}}={\bar{({\theta^1})}}$,$Z_{\bar 1}=
{\bar{({Z_1})}}$) It follows that $d{\theta}=ih_{1{\bar 1}}{\theta^1}
{\wedge}{\theta^{\bar 1}}$ for some nonzero real function $h_{1{\bar 1}}$.
If both $J$ and $\theta$ are oriented, then $h_{1{\bar 1}}$ is positive.
In this case we call such a pseudohermitian structure
$(J,{\theta})$ oriented, and we can choose a $Z_1$
(hence $\theta^1$) such that $h_{1{\bar 1}}
=1$. That is to say

\begin{equation}
d{\theta}=i{\theta^1}{\wedge}{\theta^{\bar 1}}.
\end{equation}

The pseudohermitian connection of $(J,{\theta})$ is the connection
${\nabla}^{{\psi}.h.}$ on $TM{\otimes}C$ (and extended to tensors) given by

\begin{eqnarray}
{\nabla}^{{\psi}.h.}Z_{1}={\omega_1}^{1}{\otimes}Z_{1},
{\nabla}^{{\psi}.h.}Z_{\bar 1}={\omega_{\bar 1}}^{\bar 1}{\otimes}Z_{\bar 1},
{\nabla}^{{\psi}.h.}T=0 \nonumber
\end{eqnarray}

\noindent in which the 1-form ${\omega_1}^{1}$ is uniquely determined
by the following equation with a normalization condition:

\begin{eqnarray}
&d{\theta^1}={\theta^1}{\wedge}{\omega_1}^{1}+{A^1}_{\bar 1}{\theta}
{\wedge}{\theta^{\bar 1}}& \\
& {\omega_1}^{1}+{\omega_{\bar 1}}^{\bar 1}=0. & \nonumber
\end{eqnarray}
 
The coefficient ${A^1}_{\bar 1}$ in (5.2) is called the (pseudohermitian)
torsion. Since $h_{1{\bar 1}}=1$, $A_{{\bar 1}{\bar 1}}=h_{1{\bar 1}}
{A^1}_{\bar 1}={A^1}_{\bar 1}$. And $A_{11}$ is just the complex
conjugate of $A_{{\bar 1}{\bar 1}}$. Differentiating ${\omega_1}^{1}$
gives

\begin{equation}
d{\omega_1}^{1}={\cal W}{\theta^1}{\wedge}{\theta^{\bar 1}}
+2iIm(A_{11,{\bar 1}}{\theta^1}{\wedge}{\theta})
\end{equation}

\noindent where $\cal W$ is the Tanaka-Webster curvature. Write
${\omega_1}^{1}=i{\omega}$ for some real 1-form $\omega$ by the
second condition of (5.2). This $\omega$ is just the one used in
previous sections. Write $Z_{1}=\frac{1}{2}(e_{1}-ie_{2})$ for real
vectors $e_{1},e_{2}$. Now the real version of (5.3) reads:

\begin{equation}
d{\omega}(e_{1},e_{2})=-2{\cal W}.
\end{equation}

Let $e^{1}=Re({\theta^1}),e^{2}=Im({\theta^1})$. Then $\{ e^{0}={\theta},
e^{1},e^{2} \}$ is dual to $\{ e_{0}=T,e_{1},e_{2} \}$. The
oriented pseudohermitian structure $(J,{\theta})$ induces a Riemannian
structure $h$ on $\xi$: $h(u,v)=\frac{1}{2}d{\theta}(u,Jv)$. The adapted
metric of $(J,{\theta})$ is the Riemannian metric on $TM$ defined by
${\theta}^{2}+h=(e^{0})^{2}+(e^{1})^{2}+(e^{2})^{2}$, still denoted $h$.
The Riemannian connection forms ${\omega}_{j}^{i}$ are uniquely
determined by the following equations:

\begin{eqnarray}
& de^{i}=e^{j}{\wedge}{\omega}_{j}^{i} & \\
& {\omega}_{j}^{i}+{\omega}_{i}^{j}=0. & \nonumber
\end{eqnarray}

Comparing (5.5) with (5.1),(5.2), we can relate ${\omega}_{j}^{i}$ to
the pseudohermitian connection ${\omega_1}^{1}=i{\omega}$ and torsion
${A^1}_{\bar 1}={\lambda}+i{\mu}$ (${\lambda},{\mu}$ being real) 
as follows:

\begin{eqnarray}
&&{\omega}_{1}^{2}={\omega}+{\theta},\\
&&\left\{ \begin{array}{c}
          {\omega}_{0}^{1}={\lambda}e^{1}+({\mu}-1)e^{2}\\
          {\omega}_{0}^{2}=({\mu}+1)e^{1}-{\lambda}e^{2}.
          \end{array} \right.
\end{eqnarray}
          
Observe that (5.1) and (5.2) imply

\begin{eqnarray}
&&(\frac{-i}{2})[e_{1},e_{2}]=[Z_{\bar 1},Z_{1}]=iT+{\omega_1}^{1}
(Z_{\bar 1})Z_{1}-{\omega_{\bar 1}}^{\bar 1}(Z_{1})Z_{\bar 1},\\
&&[Z_{\bar 1},T]={A^1}_{\bar 1}Z_{1}-{\omega_{\bar 1}}^{\bar 1}(T)
Z_{\bar 1}.
\end{eqnarray}

\bigskip
\begin{tabular}{ll}
\hbox{Jih-Hsin Cheng} & \hbox{Hung-Lin Chiu} \\
\hbox{Institute of Mathematics} & \hbox{Center for General Education} \\
\hbox{Academia Sinica, Taipei, R.O.C.} & \hbox{Dahan College of Engineering and Business} \\
\hbox{email: cheng@math.sinica.edu.tw} & \hbox{Hualien, Taiwan, R.O.C.}\\
 & \hbox{email: hlchiu@mss.dahan.edu.tw}
\end{tabular}
}

\end{document}